\documentstyle[12pt]{article}
\begin{document}
\begin{flushright}
Mannheimer Manuskripte, Nr. 240
\end{flushright}
\vspace{0.7cm}
\begin{center}
{\Large On Fedosov's approach to Deformation Quantization with
Separation of Variables}\\
{\bf Alexander V. Karabegov}
\footnote{the author was supported by the Alexander von Humboldt
Foundation}\\
Department of Mathematics and Computer Science\\
University of Mannheim, D7, 27\\
D-68131 Mannheim, Germany\\
E-mail: kara@math.uni-mannheim.de
\footnote{on leave of absence from the Joint Institute for Nuclear
Research, Dubna 141980, Moscow region, Russia.
E-mail: karabegov@vxjinr.jinr.ru}
\end{center}
\begin{abstract}

The description of all deformation quantizations with separation of
variables on a K\"ahler manifold from [8] is used to identify the Fedosov
star-product of Wick type constructed by M. Bordemann and S. Waldmann in
[3]. This star-product is shown to be the one with separation of
variables which corresponds to the trivial deformation of the K\"ahler
form in the sense of [8].  To this end a formal Fock bundle on a K\"ahler
manifold is introduced and an associative multiplication on its sections
is defined. 

\end{abstract}

{\bf Introduction}

For a given vector space $E$ we call formal vectors the elements of the
space $E[\nu^{-1},\nu]]$ of formal Laurent series in a formal parameter
$\nu$ with a finite principle part and coefficients in $E$.  Thus we
consider the field of formal numbers ${\bf K}={\bf C}[\nu^{-1},\nu]]$,
formal functions, forms and differential operators. 

Deformation quantization of a Poisson manifold $(M,\{\cdot,\cdot\})$, as
defined by Bayen, Flato, Fronsdal, Lichnerowicz and Sternheimer in
[2], is a structure of associative algebra on the space of formal
functions ${\cal F}=C^\infty(M)[\nu^{-1},\nu]]$. The product $\ast$ in
this algebra (called a star-product) is a ${\bf K}$-linear $\nu$-adically
continuous product given on functions $f,g\in C^\infty(M)$ by the formula

\begin{equation}
f\ast g=\sum_{r=0}^\infty \nu^r C_r(f,g).
\end{equation}

In (1) $C_r$ are bidifferential operators such that $C_0(f,g)=fg,\
C_1(f,g)-C_1(g,f)=i\{f,g\}$. The constant $1$ is assumed to be the unit in
the algebra $({\cal F},\ast)$.

Two star-products $\ast_1$ and $\ast_2$ are called equivalent if there
exists an isomorphism of algebras $B:({\cal F},\ast_1)\to ({\cal
F},\ast_2)$ given by a formal differential operator $B=1+\nu B_1+\nu^2
B_2+\dots$.

The problem of existence and classification up to equivalence of
star-products on Poisson manifolds was first solved for
symplectic manifolds (the main references are
[5,6,7,12,13]; for a historical account see [14]). In
the general case it was solved by Kontsevich [10].

Let $M$ be a K\"ahler manifold, endowed with a K\"ahler $(1,1)$-form
$\omega_{-1}$ and the corresponding Poisson bracket.  In [8] we gave a
simple geometric description of all star-products on $M$ which have the
following property of separation of variables: in a local holomorphic
chart the operators $C_r$ from (1)  act on the first argument by
antiholomorphic derivatives, and on the second argument by holomorphic
ones. We have shown that these star-products are naturally parametrized
by geometric objects, the formal deformations of the K\"ahler form
$(1/\nu)\omega_{-1}$. 

The interest in deformation quantization with separation of
variables is explained by the fact that the Wick star-product on
${\bf C}^n$ and the star-products obtained from Berezin's
quantization on K\"ahler manifolds in [4,11,9]
have the property of separation of variables.

In [3] Bordemann and Waldmann constructed a star-product with separation
of variables on an arbitrary K\"ahler manifold $(M,\omega_{-1})$, using
the geometric approach developed by Fedosov in [6,7]. The goal of this
letter is to identify the star-product obtained in [3], using the
parametrization from [8]. We show that this star-product corresponds to
the trivial deformation of the K\"ahler form $(1/\nu)\omega_{-1}$.

{\bf 1. Deformation quantizations with separation of variables}

For an open subset $U\subset M$ set ${\cal
F}(U)=C^\infty(U)[\nu^{-1},\nu]]$. Since the star-product (1) is given
by formal bidifferential operators, it can be localized to any open subset
$U\subset M$. We denote its restriction to ${\cal F}(U)$ also by $\ast$. 

Denote by ${\cal L}^\ast (U)$ and ${\cal R}^\ast (U)$ the sets of all
operators of left and right star-multiplica\-ti\-on in the algebra $({\cal
F}(U),\ast)$ respectively. All these operators are formal differential
ones. The subalgebras ${\cal L}^\ast (U)$ and ${\cal R}^\ast (U)$ of the
algebra of formal differential operators on $U$ are commutants of each
other.

Now let $(M, \omega_{-1})$ be a K\"ahler manifold with the K\"ahler
$(1,1)$-form $\omega_{-1}$.  Consider a star-product $\ast$ on $M$ with
the following property of separation of variables.  For an arbitrary
local coordinate chart $U\subset M$ with holomorphic coordinates
$\{z^k\}$ (and antiholomorphic coordinates $\{\bar z^l\}$)  assume that
the operators from ${\cal L}^\ast (U)$ contain only holomorphic
derivatives and the operators from ${\cal R}^\ast (U)$ contain only
antiholomorphic ones. This is equivalent to the fact that the operators
from ${\cal L}^\ast (U)$ and ${\cal R}^\ast (U)$ commute with the
point-wise multiplication operators by antiholomorphic and holomorphic
functions on $U$ respectively. It means that, given a holomorphic
function $a$ and antiholomorphic function $b$ on $U$, the point-wise
multiplication operators by $a$ and $b$ belong to ${\cal L}^\ast (U)$ and
${\cal R}^\ast (U)$ respectively. Therefore $L_a^\ast =a$ and $R_b^\ast
=b$, so that for $f\in{\cal F}(U)\ a\ast f=af,\ f\ast b=bf$ holds.  This
property was used for the definition of quantization with separation of
variables in [8]. 

It was shown in [8] that the star-products with separation of variables
on $(M,\omega_{-1})$ are in 1---1 correspondence with the formal
deformations of the K\"ahler form $(1/\nu)\omega_{-1}$, i.e., with the
formal forms $\omega=(1/\nu)\omega_{-1}+\omega_0+\nu\omega_1+\dots$ such
that all $\omega_r,\ r\geq 0$, are closed but not necessarily
nondegenerate $(1,1)$-forms on $M$.  

Given an arbitrary formal deformation $\omega$ of the K\"ahler form
$(1/\nu)\omega_{-1}$, one can recover the corresponding star-product with
separation of variables as follows. On each contractible coordinate chart
$(U,\{z^k\})$ on $M$ choose a formal potential
$\Phi=(1/\nu)\Phi_{-1}+\Phi_0+\nu\Phi_1+\dots$ of the form $\omega$, so
that $\omega=i\partial\bar\partial\Phi$.  Then $L_{\partial\Phi/\partial
z^k}^\ast = \partial\Phi/\partial z^k+\partial/\partial z^k$ and
$R_{\partial\Phi/\partial\bar z^l}^\ast = \partial\Phi/\partial\bar z^l+
\partial/\partial\bar z^l$.  Moreover, the set ${\cal L}^\ast (U)$
consists of all formal differential operators which commute with all
$R_{\bar z^l}^\ast =\bar z^l$ and $R_{\partial\Phi/\partial\bar z^l}^\ast
= \partial\Phi/\partial\bar z^l+ \partial/\partial\bar z^l$, and,
respectively, ${\cal R}^\ast (U)$ is the commutant of the set of all
operators $L_{z^k}^\ast =z^k$ and $L_{\partial\Phi/\partial z^k}^\ast$. 
This completely determines the star-product.

{\it Remark.} In [3] star-products with separation of variables on
K\"ahler manifolds are called star-products of Wick type, since the Wick
star-product is the simplest one of this kind. However, one can consider
star-products with separation of variables on an arbitrary symplectic
manifold endowed with a pair of transversal Lagrangean polarizations (see
[1]). In the K\"ahler case these are the holomorphic and antiholomorphic
polarizations.

{\bf 2. The formal Wick algebras bundle and the formal Fock bundle}

Consider ${\bf C}^n$ with holomorphic coordinates $\{\zeta^k\}$ (and
antiholomorphic coordinates $\{\bar\zeta^l\}$)  endowed with a Hermitian
(1,1)-form $ig_{kl}d\zeta^k\land d\bar\zeta^l$ (here $g_{kl}$ are
constants). Denote by $\circ$ the Wick star-product on $({\bf
C}^n,ig_{kl}d\zeta^k\land d\bar\zeta^l)$. This is the star-product with
separation of variables, corresponding to the trivial deformation of the
(1,1)-form $(1/\nu)ig_{kl}d\zeta^k\land d\bar\zeta^l$. The Wick
star-product of functions $f,g\in C^\infty({\bf C}^n)$ is given by the
well-known explicit formula $$
    f\circ g=\sum_{r=0}^\infty \frac{\nu^r}{r!}g^{l_1k_1}\dots g^{l_rk_r}
\frac{\partial^r f}{\partial\bar\zeta^{l_1}\dots\partial\bar\zeta^{l_r}}
\frac{\partial^r g}{\partial\zeta^{k_1}\dots\partial\zeta^{k_r}},
$$ 
where $(g^{lk})$ is the matrix inverse to $(g_{kl})$.
Here, as well as in the rest of the letter we use Einstein's summation
convention.

Introduce the following gradings on the variables
$\nu,\zeta^k,\bar\zeta^l:\ 
deg_\nu(\nu)=1,\ \deg_\nu(\zeta)=\deg_\nu(\bar\zeta)=0;\
deg'_s(\zeta)=1,\ deg'_s(\nu)=deg'_s(\bar\zeta)=0;\
deg''_s(\bar\zeta)=1,\ deg''_s(\nu)=deg''_s(\zeta)=0;\
deg_s=deg'_s+deg''_s;\
Deg'=deg_\nu+deg'_s;\
Deg''=deg_\nu+deg''_s;\
Deg=Deg'+Deg''=2 deg_\nu+deg_s.$

The Wick product $\circ$ is a graded product on polynomials in
$\nu,\zeta^k,\bar\zeta^l$ with respect to the gradings $Deg',\ Deg''$ and
$Deg$. The total grading $Deg$ is analogous to the one on the formal Weyl
algebra used by Fedosov. 

The "normal ordering" procedure establishes a 1---1 correspondence
between the polynomials from ${\bf K}[\zeta^k,\bar\zeta^l]$ and
holomorphic differential operators on ${\bf C}^n$ with coefficients in
${\bf K}[\zeta^k]$. Set $\hat\zeta^k=\zeta^k,\ \hat{\bar\zeta^l}=\nu
g^{lk}\partial/\partial\zeta^k$. The "normal ordering" relates to a
polynomial
$\phi(\zeta,\bar\zeta)=\phi_{\alpha,\beta}\zeta^\alpha\bar\zeta^\beta$
the operator
$\hat\phi=\phi_{\alpha,\beta}\hat\zeta^\alpha\hat{\bar\zeta}^\beta$. 
Here $\alpha=(k_1,\dots,k_p),\beta=(l_1,\dots,l_q)$ are multi-indices,
$\zeta^\alpha=\zeta^{k_1}\dots\zeta^{k_p},\
\bar\zeta^\beta=\bar\zeta^{l_1}\dots\bar\zeta^{l_q},\
\hat\zeta^\alpha=\hat\zeta^{k_1}\dots\hat\zeta^{k_p},
\hat{\bar\zeta}^\beta=\hat{\bar\zeta}^{l_1}\dots\hat{\bar\zeta}^{l_q}$
and $\phi_{\alpha,\beta}\in{\bf K}$ is symmetric with respect to $\alpha$ 
and $\beta$ separately. The polynomial $\phi$ is called the
Wick symbol of the operator $\hat\phi$. The operator product transferred
to Wick symbols provides the Wick product $\circ$. 

The Wick product $\circ$ can be extended to the space $W$ of formal series
in $\nu^{-1},\nu,\zeta^k,\bar\zeta^l$ with a finite principal part in
$\nu$,
$$
w=\sum_{r\geq r_0,p,q\geq 0}\nu^r \sum_{\alpha,\beta,
|\alpha|=p,|\beta|=q} w_{r,\alpha,\beta}\zeta^\alpha\bar\zeta^\beta.
$$ 
Here $r_0\in{\bf Z},\ \alpha=(k_1,\dots,k_p),\beta=(l_1,\dots,l_q)$ are
multi-indices, $\zeta^\alpha=\zeta^{k_1}\dots\zeta^{k_p},\
\bar\zeta^\beta=\bar\zeta^{l_1}\dots\bar\zeta^{l_q}$, and the terms of the
series are ordered by increasing degrees $Deg=p+q+2r$. Thus obtained
algebra $(W,\circ)$ is called a {\it formal Wick algebra}. 

A {\it formal Fock space} $V$ on ${\bf C}^n$ is the subspace of $W$ of
formal series in $\nu$ and $\zeta^k$, i.e., of the formal series
$v=\sum_{r\geq r_0,\alpha} \nu^r v_{r,\alpha}\zeta^\alpha$. Denote by
$\bar V$ the subspace of $W$ of formal series in $\nu$ and $\bar\zeta^l$. 

Consider the following projection operators in $W,\
\Pi'w=w|_{\bar\zeta=0},\
\Pi''w=w|_{\zeta=0}$
and $\Pi w=w|_{\zeta=\bar\zeta=0},\ w\in W$. Then $\Pi'W=V,\Pi''W=\bar V$
and $\Pi W={\bf K}$.

The kernels of the projections $\Pi'$ and $\Pi''$ consist of the formal
series $w\in W$ with all the terms containing at least one
antiholomorphic variable $\bar\zeta^l$ or a holomorphic variable
$\zeta^k$ respectively. It is easy to check that $Ker\ \Pi'$ and $Ker\
\Pi''$ are a left and a right ideals in the Wick algebra $(W,\circ)$
respectively. It follows, in particular, that $Ran\ \Pi'=V\cong W/Ker\
\Pi'$ is a left $W$-module. An element $w\in W$ acts on $V$ by a formal
holomorphic differential operator $T_w$ on ${\bf C}^n$ given by the
formula $T_wv=\Pi'(w\circ v),\ v\in V$. One can show that if $w\in{\bf
K}[\zeta,\bar\zeta]$ then $T_w=\hat w$, i.e., $T_w$ is the differential
operator with the Wick symbol $w$.  We shall say for general $w\in W$
that $w$ is the Wick symbol of $T_w$ and denote $T_w=\hat w$. It is easy
to check that the mapping $W\ni w\mapsto \hat w$ is an injective
homomorphism of the algebra $(W,\circ)$ to the algebra of formal
differential operators on ${\bf C}^n$. 

{\bf Lemma 1.} {\it For $w\in W\ \Pi'w=0$ iff the operator $\hat w$
annihilates the subspace of formal constants ${\bf K}\subset V$, and
$\Pi''w=0$ iff $Ran\ \hat w\subset Ker\ \Pi$.}

The proof of the lemma follows from elementary properties of Wick symbols.

Given a K\"ahler manifold $(M,\omega_{-1})$ of the complex dimension
$dim_{\bf C}M=n$, consider the unions of the formal Wick algebras and of
the formal Fock spaces associated to each tangent space to $M$. Thus we
obtain the bundles of formal Wick algebras ${\bf W}$ and of formal Fock
spaces ${\bf V}$ on $M$. For an open subset $U\subset M$ denote by ${\cal
W}(U)$ and ${\cal V}(U)$ the spaces of local sections of ${\bf W}$ and
${\bf V}$ on $U$ respectively. Set ${\cal W}={\cal W}(M),\ {\cal V}={\cal
V}(M)$. 

On a coordinate chart $(U,\{z^k\})$ on $M$ introduce the following
gradings on 1-forms $dz^k,d\bar z^l:\ deg'_a(dz)=deg''_a(d\bar z)=1,\
deg'_a(d\bar z)=deg''_a(dz)=0;\ deg_a=deg'a+deg''_a$.  Denote
$\Lambda=\oplus_r\Lambda^r$ the $deg_a$-graded algebra of differential
forms on $M$. 

There exist natural inclusions of the spaces ${\cal
F}\otimes\Lambda\subset{\cal V}\otimes\Lambda\subset{\cal
W}\otimes\Lambda$ of the (formal) scalar,
${\bf V}$- and ${\bf W}$-valued differential forms
on $M$ respectively (the tensor product is taken over $C^\infty(M),\
\otimes=\otimes_{C^\infty(M)}$). 

The fibrewise Wick product and the action of $W$ on $V$ in the first
factor of the tensor product together with the wedge-product of
differential forms in the second factor define the structures of
$deg_a$-graded algebra on ${\cal W}\otimes\Lambda$ and of its
$deg_a$-graded module on ${\cal V}\otimes\Lambda$.  The product in ${\cal
W}\otimes\Lambda$ will be denoted also $\circ$. The projections
$\Pi,\Pi'$ and $\Pi''$ define fibrewise projections in ${\cal
W}\otimes\Lambda$ denoted by the same symbols. The action of an element
$w\in{\cal W}\otimes\Lambda$ on the space ${\cal V}\otimes\Lambda$ is
given by the operator $\hat w$ defined, as above, by the expression $\hat
wv=\Pi'(w\circ v)$, where $v\in{\cal V}\otimes\Lambda$. We have
$\Pi'({\cal W}\otimes\Lambda)={\cal V}\otimes\Lambda$ and $\Pi({\cal
W}\otimes\Lambda)={\cal F}\otimes\Lambda$.

In the sequel we shall always denote by $\zeta^k,\bar\zeta^l$ the fiber
coordinates on the tangent bundle $TM$ in the frame $\{\partial/\partial
z^k,\partial/\partial\bar z^l\}$ on a coordinate chart $(U,\{z^k\})$ on
$M$. 

Notice that for a local section $w(z,\bar z)=\sum_{r\geq r_0,\alpha,\beta}
\nu^r w_{r,\alpha,\beta}(z,\bar z)\zeta^\alpha\bar\zeta^\beta\in{\cal
W}(U)$ the coefficients $w_{r,\alpha,\beta}(z,\bar z)$ are 
covariant tensor fields on $M$, symmetric with respect to $\alpha$ and
$\beta$ separately.

{\bf 3. Fedosov star-product of Wick type}

Recall the construction by Bordemann and Waldmann of the Fedosov
star-product of Wick type on a K\"ahler manifold $(M, \omega_{-1})$ from
[3]. (We use, however, different conventions and notations.) 

Let $\nabla$ denote the standard K\"ahler connection on $M$. It can be
naturally extended to tensors, and thus to the bundles ${\bf W}$ and ${\bf
V}$. For technical reasons it will be convenient to denote its extension
to ${\bf W}$ also by $\nabla$, and its extension to ${\bf V}$ by
$\hat\nabla$.

Express the K\"ahler form $\omega_{-1}$ on $M$ and the K\"ahler connection
$\nabla$ on ${\cal W}\otimes\Lambda$ in local coordinates $\{z^k,\bar
z^l,\zeta^k,\bar\zeta^l\}:\ \omega_{-1}=ig_{kl}dz^k\land d\bar z^l,\
\nabla=d-\Gamma_{ki}^s\zeta^i(\partial/\partial\zeta^s)dz^k
-\Gamma_{lj}^t\bar\zeta^j(\partial/\partial\bar\zeta^t)d\bar z^l$, where
$\Gamma_{ki}^s=g^{ls}\partial g_{kl}/\partial z^i$ and
$\Gamma_{lj}^t=g^{kt}\partial g_{kl}/\partial \bar z^j$ are the Kristoffel
symbols and $(g^{lk})$ is the matrix inverse to $(g_{kl})$.  Then
$\hat\nabla=d-\Gamma_{ki}^s\zeta^i(\partial/\partial\zeta^s)dz^k$. 

Introduce an element $R\in{\cal W}\otimes\Lambda^2$ such that it is given
in local coordinates $\{z^k,\bar z^l,\zeta^k,\bar\zeta^l\}$ by the
formula $R=(-g^{ts}\partial g_{kt}\land \bar\partial
g_{sl}+\partial\bar\partial g_{kl})\zeta^k\bar\zeta^l$. 

The curvature of the connection $\nabla$ on the bundle ${\bf W}$ was
calculated in [3]: $\nabla^2=(1/\nu)ad_{Wick}(R)$. A straightforward
calculation leads to the following

{\bf Lemma 2.} {\it The curvature of the connection $\hat\nabla$ on the
bundle ${\bf V}$ is expressed via $R$ as follows,
$\hat\nabla^2=(1/\nu)\hat R$. }

Introduce Fedosov's operators $\delta$ and $\delta^{-1}$ on ${\cal
W}\otimes\Lambda$. In local coordinates 
$\delta=(\partial/\partial\zeta^k)dz^k+
(\partial/\partial\bar\zeta^l)d\bar z^l$ 
and the operator $\delta^{-1}$ is defined as follows. For
an element $a\in{\cal W}\otimes\Lambda^q$ such that $deg_s=p$ set 
$\delta^{-1}a=0$ if $p+q=0$ and $\delta^{-1}a=(p+q)^{-1}(\zeta^k
i(\partial/\partial z^k)+\bar\zeta^l i(\partial/\partial \bar z^l))a$ if
$p+q>0$.

Then $\delta=(1/\nu)ad_{Wick}(\vartheta)$, where
$\vartheta=g_{kl}\bar\zeta^ldz^k-g_{kl}\zeta^kd\bar z^l$ (see [3]). 

It was shown in [3] that there exists a unique element $r\in{\cal
W}\otimes\Lambda^1$ which satisfies the equations $\delta^{-1}\ r=0$ and
$\delta r=R+\nabla r+(1/\nu)r\circ r$, and
contains only non-negative powers of $\nu$.

In [3] a flat Fedosov's connection $D$ on ${\bf W}$ is defined as
follows, $D=-\delta+\nabla+(1/\nu)ad_{Wick}(r)$. It is a $deg_a$-graded
derivation in the algebra $({\cal W}\otimes\Lambda,\circ)$. Therefore
${\cal W}_D=Ker\ D\cap{\cal W}$ is closed under Wick multiplication. 

It was proved in [3] that the mapping $\Pi:{\cal W}_D \to{\cal F}$ is, in
fact, a bijection.  Transferring the product from the Fedosov algebra
$({\cal W}_D,\circ)$ to ${\cal F}$ via this bijection, one obtains a
star-product $\ast$ on $(M,\omega_{-1})$. Moreover, it was proved in [3]
that $\ast$ is a star-product with separation of variables. The proof was
based on the following important statement (Lemma 4.5 in [3]): $r\in Ker\
\Pi'\cap Ker\ \Pi''$, i.e., in any local expression of $r$ each term
contains variables $\zeta^k$ and $\bar\zeta^l$ for some indices $k,l$. We
reformulate this statement using Lemma 1. 

{\bf Lemma 3.} {\it The operator $\hat r$ in ${\cal V}$ annihilates the
subspace ${\cal F}\otimes\Lambda\subset{\cal V}\otimes\Lambda$. In
particular, $\hat r1=0$. Moreover, $Ran\ \hat r\subset Ker\ \Pi$.}

We are going to show that the star-product with separation of variables
$\ast$ constructed in [3] corresponds to the trivial deformation
$\omega=(1/\nu)\omega_{-1}$ of the K\"ahler form $(1/\nu)\omega_{-1}$. 

{\bf 4. The Fock algebra}

Using the fact that $\delta=(1/\nu)ad_{Wick}(\vartheta)$, one can
express $D$ as follows, $D=\nabla+(1/\nu)ad_{Wick}(\gamma)$, where
$\gamma=-\vartheta+r$.

Introduce a connection $\hat D$ on ${\bf V}$ by the formula
$\hat D=\hat\nabla+(1/\nu)\hat\gamma$.

One can split the connections $\nabla,D,\hat D$, the operator $\delta$
and the element $r$ into the sums of their (1,0)- and (0,1)-components,
$\nabla=\nabla'+\nabla'',D=D'+D'',\hat D=\hat D'+\hat D'',
\delta=\delta'+\delta'',r=r'+r''$.  In local coordinates denote
$\nabla_k=\nabla_{\partial/\partial z^k},\
\nabla_l=\nabla_{\partial/\partial \bar z^l}$, so that
$\nabla'=\nabla_kdz^k,\ \nabla''=\nabla_ld\bar z^l$. Introduce similarly
$D_k,D_l,\hat D_k,\hat D_l$. Let $r=r_kdz^k+r_ld\bar z^l$ be a local
expression of the element $r$. Then $r'=r_kdz^k,r''=r_ld\bar z^l$.

A simple calculation shows that
$(1/\nu)\hat{\vartheta}=\partial/\partial\zeta^kdz^k-\eta_ld\bar z^l$,
where $\eta_l=(1/\nu)g_{kl}\zeta^k$. Therefore, 
\begin{equation}
\hat D_k=\partial/\partial z^k-\partial/\partial
\zeta^k+(1/\nu)\hat r_k;\
\hat D_l=\partial/\partial\bar z^l+\eta_l+(1/\nu)\hat r_l.
\end{equation}

{\bf Lemma 4.} {\it Let $f\in{\cal F}(U)$, where $(U,\{z^k\})$ is a
coordinate chart on $M$. Then $\hat D_kf=\partial f/\partial z^k$. In
particular, $\hat D_k1=0$.}

The lemma trivially follows from Lemma 3 and formula (2).

{\bf Lemma 5.} {\it For $w\in{\cal W}\otimes\Lambda$ one has
$[\hat\nabla,\hat w]=\widehat{\nabla w}$.} 

Here, as well as below, the commutator is the $deg_a$-graded commutator in
the graded algebra of endomorphisms of ${\cal V}\otimes\Lambda$. 

The lemma is an easy consequence of the fact that $\nabla$ is a
$deg_a$-graded derivation of the algebra $({\cal
W}\otimes\Lambda,\circ)$. It implies the following

{\bf Proposition 1.} {\it For $w\in{\cal W}\otimes\Lambda$ the
formula $[\hat D,\hat w]=\widehat{Dw}$ holds.}

Denote $\omega=(1/\nu)\omega_{-1}$.

{\bf Lemma 6.} {\it\\
(i) $[\hat \nabla,\hat{\vartheta}]=0$;\\
(ii) $(1/\nu)[\hat{\vartheta},\hat r]=\widehat{\delta r}$;\\
(iii) $\hat{\vartheta}^2=i\nu^2\omega$.} 

Lemma is proved by straightforward calculations. It implies the following

{\bf Proposition 2.} {\it The connection $\hat D$ on ${\bf V}$ has a
scalar curvature, $\hat D^2=i\omega$.}
 
The subspace ${\cal W}_{D''}=Ker\ D''\cap {\cal W}$ of the algebra $({\cal
W},\circ)$ is closed under the Wick product. We shall use the algebra
$({\cal W}_{D''},\circ)$ to define a product on the space ${\cal V}$.

Introduce Fedosov's operator $\delta''^{-1}$ on ${\cal W}\otimes\Lambda$
defining it in the local coordinates on a chart $(U,\{z^k\})$ as follows.
Let $w\in ({\cal W}\otimes\Lambda)(U)$ be such that $deg''_s(w)=p,\
deg''_a(w)=q$. Set $\delta''^{-1}a=0$ if $p+q=0$ and
$\delta''^{-1}a=(p+q)^{-1}\bar\zeta^l i(\partial/\partial \bar z^l)a$ if
$p+q>0$. Then for $w\in{\cal W}\otimes\Lambda$ one has
$(\delta''\delta''^{-1}+\delta''^{-1}\delta'')w=w-w_0$, where $w_0$ is
the $(deg''_s+deg''_a)$-homogeneous component of $w$ of the degree 0. 
For an element $w\in{\cal W}\otimes\Lambda$ denote by $w^{(q)}$ its
$Deg''$-homogeneous component of the degree $q$. 

The following proposition can be proved by Fedosov's technique
developed in [6].

{\bf Proposition 3.} {\it The mapping $\Pi':{\cal W}_{D''}\to{\cal V}$ is
a bijection. For an element $v\in {\cal V}$ such that $deg_\nu(v)=0$
(i.e., which does not depend on the
formal parameter $\nu$) the unique element $w\in{\cal W}_{D''}$ such that
$v=\Pi'w$ can be calculated recursively with respect to the degree $Deg''$
by 
$$
   w^{(0)}=v;
$$
$$
w^{(q+1)}=\delta''^{-1}(\nabla''w^{(q)}+
(1/\nu)\sum_{p=0}^qad_{Wick}(r''^{(p+1)})w^{(q-p)}). 
$$ }

Denote by $\bullet$ the product in ${\cal V}$ obtained by pushing forward
the product in the algebra $({\cal W}_{D''},\circ)$ by the mapping
$\Pi'$.  Thus we obtain a {\it Fock algebra} $({\cal V},\bullet)$. For
$v\in{\cal V}$ denote by $L_v^\bullet,R_v^\bullet$ the operators of left
and right multiplication by $v$ in the algebra $({\cal V},\bullet)$
respectively.  Set ${\cal L}^\bullet=\{L^\bullet_v|v\in{\cal V}\},\ {\cal
R}^\bullet=\{R^\bullet_v|v\in{\cal V}\}$. 

{\bf Lemma 7.} {\it For $w\in{\cal W}_{D''}$ the operator $\hat w$
coincides
with the left multiplication operator by the element $v=\Pi'w$ in the
Fock algebra $({\cal V},\bullet),\ \hat w=L^\bullet _v$. 

Proof.} For $w_1,w_2\in{\cal W}_{D''}$ set $v_1=\Pi'w_1,\ v_2=\Pi'w_2$. 
Then, by definition, $v_1\bullet v_2=\Pi'(w_1\circ w_2)$.  Since $\Pi'$
is a projection, $w_2-v_2\in Ker\ \Pi'$.  Taking into account that $Ker\
\Pi'$ is a left ideal in the algebra $({\cal W},\circ)$, we get
$w_1\circ(w_2-v_2)\in Ker\ \Pi'$. Therefore $\Pi'(w_1\circ
w_2)=\Pi'(w_1\circ v_2)=\hat w_1v_2$, whence the Lemma follows.  $\Box$

Since the action of the operators $\hat w,\ w\in{\cal W}$, on ${\cal V}$
is fibrewise, it follows from Lemma 7 that the operator of point-wise
multiplication by $f\in{\cal F}$ (also denoted by $f$) commutes with all
operators from ${\cal L}^\bullet$. Therefore, $f\in{\cal R}^\bullet$,
namely, $R^\bullet_f=f$. 

Fix a coordinate chart $(U,\{z^k\})$ on $M$.

{\bf Lemma 8.} $R^\bullet_{\eta_l}=\hat D_l$.

{\it Proof.} Let $w\in{\cal W}_{D''}(U),\ v=\Pi'w\in{\cal V}(U)$. It
follows from Lemma 7 and Proposition 1 that $[\hat D_l,L_v^\bullet]=[\hat
D_l,\hat w]=\widehat{D_lw}=0$, therefore $\hat D_l\in{\cal R}^\bullet$. 
Using formula (2) and Lemma 3 we get $\hat D_l1=\eta_l$, whence $\hat
D_l=R^\bullet_{\eta_l}$.  $\Box$

Denote ${\cal U}=\Pi'({\cal
W}_D)\subset{\cal V}$. Since ${\cal W}_D\subset{\cal W}_{D''}$, and the
projection $\Pi'$ establishes an isomorphism of the algebras $({\cal
W}_{D''},\circ)$ and $({\cal V},\bullet)$, the subspace ${\cal
U}\subset{\cal V}$ is closed under multiplication $\bullet$ and the
projection $\Pi'$ maps the Fedosov algebra $({\cal W}_D,\circ)$
isomorphically onto the subalgebra $({\cal U},\bullet)$ of the Fock
algebra $({\cal V},\bullet)$. 

{\bf Lemma 9.} {\it For $w\in{\cal W}_{D''}(U)$ and $v=\Pi'w\in {\cal
V}(U)$ one has $D_kw\in{\cal W}_{D''}(U)$ and $[\hat
D_k,L^\bullet_v]=\widehat{D_kw}=L^\bullet_{\hat D_kv}$.

Proof.} Using Lemma 7 and Proposition 1 we obtain $[\hat
D_k,L^\bullet_v]=[\hat D_k,\hat w]=\widehat{D_kw}$.  Since Fedosov's
connection $D$ is flat, $D^2=0$, we have $[D_k,D_l]=0$, whence
$D_lD_kw=D_kD_lw=0$, i.e., $D_kw\in{\cal W}_{D''}(U)$ and therefore
$\widehat{D_kw}\in{\cal L}^\bullet(U)$. Using Lemma 4 we get $[\hat
D_k,L^\bullet_v]1=\hat D_kv-L^\bullet_v\hat D_k1=\hat D_kv$ and thus
$\widehat{D_kw}=L^\bullet_{\hat D_kv}$, which concludes the proof. $\Box$

Denote ${\cal V}_{\hat D'}(U)=Ker\ \hat D'\cap{\cal V}(U)$ the space of
local sections of the Fock bundle ${\bf V}$ on an open subset
$U\subset M$, annihilated by $\hat D'$. Set
${\cal V}_{\hat D'}={\cal V}_{\hat D'}(M)$.

{\bf Proposition 4.}  ${\cal U}={\cal V}_{\hat D'}$.

{\it Proof.} We have to show that on any coordinate chart $(U,\{z^k\})$
on $M,\ w\in{\cal W}_{D''}(U)$ and
$v=\Pi'w\in{\cal V}(U)$ the condition $D_kw=0$ holds iff $\hat D_kv=0$.
The assertion follows immediately from the equality $L^\bullet_{\hat
D_kv}=\widehat{D_kw}$ proved in Lemma 9 and the fact that the mapping
${\cal W}\ni w\mapsto\hat w$ is injective.  $\Box$

We can obtain the star-product $\ast$ on $M$ from the algebra
$({\cal V}_{\hat D'},\bullet)=({\cal U},\bullet)$. Let $v_1,v_2\in {\cal
V}_{\hat D'},\ f_1=\Pi v_1,f_2=\Pi v_2\in{\cal F}$. Then $f_1\ast
f_2=\Pi(v_1\bullet v_2)$.

Let $\Phi_{-1}$ be a local potential of the form
$\omega_{-1}=ig_{kl}dz^k\land d\bar z^l$ on a coordinate chart
$(U,\{z^k\})$ on $M$, so that $\partial^2\Phi_{-1}/\partial
z^k\partial\bar z^l=g_{kl}$. Then $\Phi=(1/\nu)\Phi_{-1}$ is a local
potential of the form $\omega=(1/\nu)\omega_{-1}$. Set
$Q_l=\partial\Phi/\partial\bar z^l+\eta_l$. 

{\bf Proposition 5.} $Q_l\in {\cal V}_{\hat D'}(U)$. 

{\it Proof.} Using Lemma 4 we get $\hat D_k\partial\Phi/\partial\bar
z^l=\partial^2\Phi/\partial z^k\partial\bar z^l=(1/\nu)g_{kl}$. It
follows from Proposition 2 that $[\hat D_l,\hat D_k]=(1/\nu)g_{kl}$. Now,
$\hat D_k\eta_l=\hat D_k\hat D_l1=\hat D_l\hat
D_k1-(1/\nu)g_{kl}=-(1/\nu)g_{kl}$ and therefore
$\hat D'Q_l=0$. $\Box$

Since $\ast$ is known to be a star-product with separation of variables,
then $R^\ast_{\bar z^l}=\bar z^l$ holds. This can be checked
also directly.  It follows from Lemma 4 that $\hat D_k\bar z^l=0$, i.e.,
$\bar z^l\in {\cal V}_{\hat D'}(U)$.  Let $v\in{\cal V}_{\hat D'}(U)$
and $f=\Pi v\in{\cal F}(U)$.  Now $f\ast\bar z^l=\Pi(v\bullet\bar
z^l)=\Pi(v\bar z^l)=f\bar z^l$, which proves the assertion. 

In order to identify the star-product with separation of variables $\ast$
it remains to calculate $R^\ast_{\partial\Phi/\partial\bar z^l}$.
Let $v\in{\cal V}_{\hat D'}(U)$ and $f=\Pi v$ as above.
Calculate first $\Pi\hat D_lv$. Using formula (2) we get $\Pi\hat
D_lv=\Pi(\partial v/\partial\bar z^l+\eta_lv+(1/\nu)\hat r_lv)$.  Since
$\Pi\eta_l=0$, we have $\Pi(\eta_lv)=0$. Lemma 3 implies that $\Pi(\hat
r_lv)=0$. Finally we obtain that $\Pi\hat D_lv=\partial f/\partial\bar
z^l$. 

Since $\Pi Q_l=\partial\Phi/\partial\bar z^l$, we get
$f\ast\partial\Phi/\partial\bar z^l=\Pi(v\bullet Q_l)=\Pi (R^\bullet_
{Q_l}v)=\Pi((\partial\Phi/\partial\bar z^l+\hat D_l)v)=
(\partial\Phi/\partial\bar z^l+\partial/\partial\bar z^l)f$.  Therefore
$R^\ast_{\partial\Phi/\partial\bar z^l}= \partial\Phi/\partial\bar
z^l+\partial/\partial\bar z^l$. Thus we have proved the desired

{\bf Theorem.} {\it The Fedosov star-product of Wick type $\ast$
on a K\"ahler manifold $(M,\omega_{-1})$
is the star-product with separation of variables corresponding to the
trivial deformation of the form $(1/\nu)\omega_{-1}$.}

One might also try to generalize the construction by Bordemann and
Waldmann to obtain arbitrary deformation quantizations with separation of
variables. 

{\bf Acknowledgements} 

The author is deeply indebted to late Professor M. Flato for the
continuous support of his research work.

The author is very grateful to B. Fedosov, A. P. Nersessian and M. 
Schlichenmaier for stimulating discussions, to the Alexander von Humboldt
foundation for the fellowship granted, and to the Department of
Mathematics and Computer Science of the University of Mannheim for their
warm hospitality.

{\bf References}

\noindent1. Astashkevich, A.: On Karabegov's quantization on semisimple
coadjoint orbits, in {\it Adv. in Geom. and Math. Phys.}, Vol. 1, eds. J.
L.  Brylinski et al., Springer-Verlag, New-York, 1998. 

\noindent 2. Bayen, F., Flato, M., Fronsdal, C., Lichnerowicz, A. and
Sternheimer, D.:  Deformation theory and quantization I, {\it Ann. of
Phys.} {\bf 111} (1978), 61-110.

\noindent 3. Bordemann, M. and Waldmann, S.: A Fedosov star product of
Wick type for K\"ahler manifolds, {\it Lett. Math. Phys.} {\bf 41}
(1997), 243-253.

\noindent 4. Cahen, M., Gutt, S. and Rawnsley, J.: Quantization of
K\"ahler manifolds, II, {\it Trans. Am. Math. Soc.} {\bf 337} (1993),
73-98.

\noindent 5. De Wilde, M. and Lecomte, P.B.: Existence of star-products
and of formal deformations of the Poisson Lie algebra of arbitrary
symplectic manifolds, {\it Lett. Math. Phys.} {\bf 7} (1983), 487-496.

\noindent 6. Fedosov, B. V.: A simple geometrical construction of
deformation quantization, {\it J. Differential Geom.} {\bf 40} (1994),
213-238.

\noindent 7. Fedosov, B. V.: {\it Deformation Quantization and Index
Theory,} Math. Top. 9, Akademie-Verlag, Berlin, 1996.

\noindent 8. Karabegov, A.V.: Deformation quantizations with separation
of variables on a K\"ahler manifold, {\it Comm. Math. Phys.} {\bf 180}
(1996), 745-755.

\noindent 9. Karabegov A. V.: Pseudo-K\"ahler quantization on flag
manifolds, {\it Comm. Math. Phys.} {\bf 200} (1999), 355-379. 

\noindent 10. Kontsevich, M.: Deformation quantization of Poisson
manifolds, q-alg/9709040. 

\noindent 11. Moreno, C.: $\ast$-products on some K\"ahler manifolds,
{\it Lett. Math.  Phys.} {\bf 11} (1986), 361-372. 

\noindent 12. Nest, R. and Tsygan, B.: Algebraic index theorem, {\it
Comm. Math. Phys.} {\bf 172} (1995), 223-262.

\noindent 13. Nest, R. and Tsygan, B.: Algebraic index theorem for
families, {\it Adv.  in Math.} {\bf 113} (1995), 151-205. 

\noindent 14. Weinstein, A.: Deformation quantization, {\it S\'eminaire
Bourbaki, 46\'eme annee, Asterisque 227,} {\bf 789} (1995), 389-409.
\end{document}